\documentclass[12pt]{article}

\usepackage {amsfonts,amsthm}
\usepackage[cp1251]{inputenc}
\usepackage[english]{babel}
\theoremstyle{plain}
\newtheorem*{Th}{Theorem}
\newtheorem{Le}{Lemma}

\newtheorem*{Prop}{Proposition}
\def\R{{\mathbb R}}
\def\C{{\mathbb C}}

\def\Z{{\mathbb Z}}

\def\e{\varepsilon}
\def\Re{{\rm Re}}
\def\Im{{\rm Im}}

\def\l{\lambda}
\def\a{\alpha}
\def\b{\beta}
\def\g{\gamma}
\def\D{\Delta}

\def\Lin{{\rm Lin}}

\def\<{\langle}
\def\>{\rangle}
\begin{document}

\title{
    \bf
Multidimensional Version of Lagrange's Problem on Mean Motion}

\author{S.\,Favorov}

\date{}

\maketitle
\begin{abstract}
The famous mean motion problem which goes back to Lagrange as follows: to prove that
any exponential polynomial with exponents on the imaginary axis has an average speed
for the amplitude, whenever the variable moves along a horizontal line. It was
completely proved by B.\,Jessen and H.\,Tornehave in Acta Math.77, 1945. Here we give
its multidimensional version.
\end{abstract}

{\it 2000 Mathematics Subject Classification:} {\small Primary 32A05, Secondary 30B50,
32A60}

{\it Keywords:} {\small  exponential polynomial, Lagrange's conjecture, mean motion}

\bigskip
Consider an exponential polynomial
\begin{equation}\label{sum}
  f(z)=\sum_{j=1}^N
  c_je^{i\l_jz},\quad c_j\in\C,\quad\l_j\in\R.
\end{equation}
 J.L.~Lagrange \cite{L} assumes that for each fixed $y\in\R$ there exist the limits
\begin{equation}\label{m1}
  c^+(y)=\lim_{\b-\a\to\infty}{\D_{\a<x<\b}\arg^+f(x+iy)\over\b-\a}
\end{equation}
and
\begin{equation}\label{m2}
  c^-(y)=\lim_{\b-\a\to\infty}{\D_{\a<x<\b}\arg^-f(x+iy)\over\b-\a}.
\end{equation}
Here $\arg^+f(x+iy)$ and $\arg^-f(x+iy)$ are branches of $\arg f(x+iy)$, which are
continuous in $x$ on every interval without zeros of $P$ and have the jumps $-p\pi$
and $+p\pi$, respectively, at zeros of multiplicity $p$, $\D^\pm_{\a<x<\b}$ are
increments of the functions $\arg^\pm f(x+iy)$ on $(\a,\,\b)$.

 J.L.~Lagrange proves his conjecture when the absolute value
of one of the coefficients in  (\ref{sum}) is greater than
 the sum of the absolute values of the other coefficients. Moreover, if this is the case, then
\begin{equation}\label{m3}
\arg^+f(x+iy)= c^+x+O(1),\quad
  \arg^-f(x+iy)= c^-x+O(1)\quad(x\to\infty),
\end{equation}
besides, the mean motions $c^+(y)$ and $c^-(y)$ are equal. He also notes that in the
case $N=2$ with arbitrary terms in (\ref{sum}), the equalities (\ref{m3}) hold as
well, but $c^+$ and $c^-$ may be different (for example, in the case $f(z)=\sin z$).

Note that (\ref{m3}) is false  for sums (\ref{sum}) in the case
$N>2$ (see F.\,Bernstein \cite{Ber}). On the other hand, H.\,Bohr
in \cite{B} proves (\ref{m3}) with $c^+=c^-$ for almost periodic
functions $f$ on $\R$ under the condition
 \begin{equation}\label{mor}
 |f|\ge\kappa>0.
 \end{equation}
 Moreover, he shows that in this case
the term $O(1)$ in (\ref{m3}) is  an almost periodic function as
well. Then B.\,Jessen \cite{J} proves that for  holomorphic almost
periodic functions in a strip $\{z=x+iy:\,\,a<y<b\}$ limits
(\ref{m1}) and (\ref{m2}) exist for all $y\in(a,\,b)$ outside of a
countable set. Moreover, he establishes a connection of  mean
motions with distributions of zeros of $f$.

 Lagrange's Conjecture for exponential polynomials is proved by
H.\,Weil \cite{W2} in the case of linearly independent $\l_1,\dots,\l_N$ over $\Z$,
and by B.\,Jessen and H.\,Tornehave \cite{JT} in the general case (see the detailed
history in the introduction to \cite{JT}). Also, they prove that the equality
$c^+(y)=c^-(y)$ for exponential polynomials remains true for all $y\in\R$ outside of a
discrete set without finite limit points.

Bohr's result holds true in the multidimensional case: for any almost periodic
function in $\R^p$ with condition (\ref{mor}) we have
 $$
 f(x)=\exp\{i\<c,x\>+g(x)\},
 $$
where $g$ is almost periodic in $x\in\R^p$, $c$ is a vector from
$\R^p$ (L.Ronkin \cite{R2}). In \cite{R1}, \cite{F}, and \cite{FG}
one can find various relations between mean motions and
distribution of zero sets for holomorphic almost periodic
functions in tube domains.

Nevertheless, I don't know papers about Lagrange's conjecture for exponential
polynomials in several variables, although there are a lot of papers devoted to
properties of such polynomials (see, for example, \cite{Ge}, \cite{K1}, \cite{K2},
\cite{P}). The result of the present paper solves Lagrange's problem in arbitrary
dimension.
\bigskip

{\bf Notations}. For $x=(x_1,\dots,x_p)\in\R^p$ put
$'x=(x_2,\dots,x_p)$, for $z=(z_1,\dots,z_p)\in\C^p$ put $z=x+iy$,
$x,y\in\R^p$. By $\< x,\,y \>$ or $\< z,\,w\>$ denote the scalar
product for vectors from $\R^p$ or the Hermitian scalar product
for $z,w \in\C^p$, $|a|$ means the Euclidean norm of the vector
$a$ in $\C^d$ or in $\R^d$. Also, by $m_d$ denote the Lebesgue
measure in $\R^d,\,d=1,2,\dots$.
\medskip

In the sequel we will use the following well known result\footnote{It can be proved
easily by induction in $d$. Also, see \cite{R3}, Ch.2, \S 2.}.
\begin{Prop}
 For any analytic set $M$ in $\Omega\subset\C^d$ and any hyperplane
$L=\{z=x+iy\in\C^d:\,y\equiv y^0\}$ we have $m_d(M\cap L)=0$.
\end{Prop}
\medskip

\begin{Th}
  For each exponential polynomial
  $$
P(z)=\sum_{j=1}^{S}
  c_je^{i\<z,\,\l^j\>},\quad c_j\in\C,\,
  \l^j\in\R^p,\,\l^j\neq\l^{j'}\,\hbox{for}\,j\neq j'.
   $$
and each $y\in\R^p$, there exist the limits
 $$
\lim_{\min_j(\b_j-\a_j)\to\infty} \prod_{1\le j\le p}
(\b_j-\a_j)^{-1}
\int_{\Pi^{(p-1)}(\,'\a,\,'\b)}\D_{\a_1<x_1<\b_1}\arg^+
P(x+iy)\,d\;'x
 $$
and
 $$
\lim_{\min_j(\b_j-\a_j)\to\infty} \prod_{1\le j\le p}
(\b_j-\a_j)^{-1}
\int_{\Pi^{(p-1)}(\,'\a,\,'\b)}\D_{\a_1<x_1<\b_1}\arg^-
P(x+iy)\,d\;'x.
 $$
Here $\a,\,\b\in\R^p,\,\a_j<\b_j \,\forall j$, and
$\Pi^{(p-1)}('\a,'\b)=\{'x\in\R^{p-1}:\,\a_j<x_j<\b_j,\,
j=2,\dots,p\}.$
\end{Th}

 Similar to the one-dimensional case, the proof of the Theorem is based
 on the following lemma.
\begin{Le}[for $p=1$ see \cite{W1}]\label{1}
 Suppose $g(u),\, u=(u_1,\dots,u_N),$ is a $2\pi$-periodic function
in each variable $u_1,\dots,u_N$,  and  $\mu^1,\dots,\mu^N\in\R^p$ are linearly
independent over $\Z$. If $g(u)$ is integrable in the sense of Riemann on
$[0,\,2\pi]^N$, then the limit
\begin{equation}\label{m4}
\lim_{\min_j(\b_j-\a_j)\to\infty}\prod_j(\b_j-\a_j)^{-1}
\int_{\Pi^{(p)}(\a,\b)}g(\<\mu^1,x\>,\dots,\<\mu^N,x\>)\,dx
\end{equation}
exists and equals the average {\bf M}$g$ of the function $g(u)$
over the cube $[0,\,2\pi]^N$.
\end{Le}

 \noindent{\bf Proof.} The proof of the  Lemma is the same
as in the one--dimensional case.

We may assume that $g$ is a real-valued function.
 If $g$ is a trigonometric polynomial of the form
 \begin{equation}\label{pol}
\sum_{k\in\Z^N} b_{k} e^{i\<k,\,u\>},
 \end{equation}
then its average equals the coefficient $b_0$. Since $k_1
\mu^1+\dots+k_N\mu^N=0$ only for the case $k=(k_1,\dots,k_N)=0$,
we see that (\ref{m4}) equals $b_0$ as well.

Furthermore, an arbitrary continuous $2\pi$-periodic in each
variable function $g$ can be uniformly approximated by polynomials
(\ref{pol}), hence we obtain the conclusion of the Lemma in this
case too.

Finally, for any integrable in the sense of Riemann function $g$
and any $\e>0$ there are continuous $2\pi$-periodic in each
variable functions $g_\e(u)\le g(u)$ and $g^\e(u)\ge g(u)$ such
that
 $$
{\bf M}g^e\le {\bf M}g+\e,\quad {\bf M}g_e\ge {\bf M}g-\e.
 $$
 Then we get
 $$
 \limsup_{\min_j(\b_j-\a_j)\to\infty}\prod_j(\b_j-\a_j)^{-1} \int_{\Pi^{(p)}(\a,\b)}
 g(\<\mu^1,x\>,\dots,\<\mu^N,x\>)\,dm_p(x)\le{\bf M}g^\e,
 $$
 $$
 \liminf_{\min_j(\b_j-\a_j)\to\infty}\prod_j(\b_j-\a_j)^{-1} \int_{\Pi^{(p)}(\a,\b)}
 g(\<\mu^1,x\>,\dots,\<\mu^N,x\>)\,dm_p(x)\ge{\bf M}g_\e.
 $$
 Since $\e$ is arbitrary small, we obtain the assertion of the
 Lemma.
\medskip

We also need the following simple assertion.
\begin{Le}\label{2}
For any real numbers $\g_1,\dots,\g_n$ there is a constant
$C<\infty$ such that the number of zeros in the segment $[-1,\,1]$
of an arbitrary exponential polynomial $g(s)\not\equiv0$ of the
form
 $$
q(s)=\sum_{k=1}^n a_ke^{i\g_ks},\quad a_k\in\C,
 $$
 does not exceed $C$.
\end{Le}
{\bf Proof of the Lemma}. Collecting similar terms, rewrite $q(s)$
in the form
 $$
q(s)=\sum_{k=1}^{n'} a'_ke^{i\g'_ks},\quad
\g'_k\neq\g'_l\quad\hbox {for}\quad k\neq l.
 $$
We may suppose $\sum_k |a'_k|=1$. The functions $e^{i\g'_ks}$ are linearly independent
over $\C$. Hence, if at least one of the coefficient $a'_k$ does not vanish, then
$g(s)\not\equiv0$.  Using Hurwitz' theorem, we obtain an easy proof of the Lemma by
contradiction.
\bigskip

 {\bf Proof of the Theorem}.

 Let $\mu^1,\dots,\mu^N$ be a basis of
the group $\Lin_\Z \{\l^1,\dots,\l^S\}$. Therefore,
 $$
\l^j=\sum_{r=1}^N k_{r,j}\mu^r,\quad k_{r,j}\in\Z,\quad 1\le j\le
S,\,1\le r\le N.
 $$
Set for $w=(w_1,\dots,w_N)\in\C^N$
 $$
 F(z,w) =\sum_{j=1}^S
  c_j\exp\{i\<z,\l^j\>+i\sum_{r=1}^N k_{r,j}w_r\}.
 $$
The function $F(z,\,u),\,u\in\R^N$, is $2\pi$-periodic in each
variable $u_1,\dots,u_N$ and
 \begin{equation}\label{shift}
F(t+iy,\<\mu^1,x\>,\dots,\<\mu^N,x\>)=P(x+iy+t)\quad\forall\,t\in\R^p.
 \end{equation}

Fix $y=y^{(0)}\in\R^p$.  Note that the set
$$
M=\{w\in\C^N:\,F(z_1,i'y^{(0)},w)=0 \quad\forall\,z_1\in\C\}
$$
$$
=\{w\in\C^N:\,0=F(0,i'y^{(0)},w)=F'_{z_1}(0,i'y^{(0)},w)=
F''_{z_1}(0,i'y^{(0)},w)=\dots\}
$$
is closed and analytic in $\C^N$, therefore, by the Proposition,
$m_N(M\cap\R^N)=0$. Consider the functions
 $$
I^+(u)=\D_{-1/2<x_1<1/2}\arg^+ F(x_1+iy_1^{(0)},i'y^{(0)},u).
 $$
 and
  $$
I^-(u)=\D_{-1/2<x_1<1/2}\arg^- F(x_1+iy_1^{(0)},i'y^{(0)},u).
 $$
We will suppose that $I^+(u)\equiv I^-(u)\equiv0$ if $u\in
M\cap\R^N$.

Let us prove that these functions are uniformly bounded and
continuous almost everywhere in $u\in[0,\, 2\pi]^N$.

If $F(x_1+iy_1^{(0)},i'y^{(0)},u^{(0)})\neq0$ for all
$x_1\in[-1/2,\,1/2]$, then  the function
$F_z'(x_1+iy_1^{(0)},i'y^{(0)},u)/F(x_1+iy_1^{(0)},i'y^{(0)},u)$
is continuous and uniformly bounded in $x\in[-1/2,\,1/2]$ and $u$
belonging to a neighborhood of $u^{(0)}$. Hence the functions
 $I^+(u)$ and $I^-(u)$ are equal, uniformly bounded and continuous
in this neighborhood.

Suppose that $u^{(0)}\not\in M\cap\R^N$ and
$F(x_1^{(1)}+iy_1^{(0)},i'y^{(0)},u^{(0)})=0$ at a point $x_1^{(1)}\in[-1/2,\,1/2]$.
Since $F(z_1,i'y^{(0)}, u^{(0)})\not\equiv 0$, we can use the Weierstrass Preparation
Theorem  (see, for example, \cite{H}). Hence there are $\e>0$, $\delta>0$, and
pseudopolynomial
 \begin{equation}\label{pol1}
 P_1(z_1,i'y^{(0)},w)= (z_1-x_1^{(1)}-iy_1^{(0)})^r + a_1(w)(z_1-x_1^{(1)}-iy_1^{(0)})^{r-1} +\dots+a_r(w)
 \end{equation}
with holomorphic in the ball $\{w:\, |w-u^{(0)}|<\e\}$
coefficients $a_j(w)$ such that
\begin{equation}\label{co}
 a_j(u^{(0)})=0,\quad j=1,\dots,r,
 \end{equation}
and
 $$
   F(z_1,i'y^{(0)},w)=P_1(z_1,i'y^{(0)},w)F_1(z_1,i'y^{(0)},w),\quad
   F_1(z_1,i'y^{(0)},w)\neq 0
 $$
in the set $\{(z_1,w):\,|w-u^{(0)}|<\e,\,
|z_1-x_1^{(1)}-iy_1^{(0)}|<\delta\}$.

Note that for a sufficiently small $\e$ each zero of
$P_1(z_1,i'y^{(0)},w)$ belongs to the disc
$|z_1-x_1^{(1)}-iy_1^{(0)}|<\delta$. Hence the function
$F_1=F/P_1$ is holomorphic in the set
$\{(z_1,w):\,z_1\in\C,\,|w-u^{(0)}|<\e\}$.

Let $x_1^{(2)}$ be another point of the segment $[-1/2,1/2]$ such
that $F(x_1^{(2)}+iy_1^{(0)},i'y^{(0)},u^{(0)})=0$. Using the
Weierstrass Preparation Theorem for $F_1(z_1,i'y^{(0)},w)$ in a
neighborhood of the point $(x_1^{(2)}+iy_1^{(0)},u^{(0)})$, we get
 $$
F_1(z_1,i'y^{(0)},w)=P_2(z_1,i'y^{(0)},w)F_2(z_1,i'y^{(0)},w).
 $$
Here $P_2(z_1,i'y^{(0)},w)$ has the form (\ref{pol1}) with
$x_1^{(2)}$ instead of $x_1^{(1)}$,  the function
$F_2(z_1,i'y^{(0)},w)$ is holomorphic in the set
$\{(z_1,w):\,z_1\in\C,\,|w-u^{(0)}|<\e\}$ and has not zeros in the
set $\{(z_1,w):\,|w-u^{(0)}|<\e,\,
|z_1-x_1^{(2)}-iy_1^{(0)}|<\delta\}$. Continuing in the same way,
we get the representation
 \begin{equation}\label{rep}
    F(z_1,i'y^{(0)},w)=P_1(z_1,i'y^{(0)},w)\cdot\dots\cdot
    P_s(z_1,i'y^{(0)},w)G(z_1,i'y^{(0)},w),
 \end{equation}
where the pseudopolynomials $P_j$ have form (\ref{pol1}) with various points $\tilde
x_1\in[-1/2,1/2]$ instead of  $x_1^{(1)}$, their coefficients satisfy (\ref{co}), and
the holomorphic function $G(z_1,i'y^{(0)},w)$ does not vanish in a neighborhood of the
set $\{(z_1,w):\,x_1\in[-1/2,1/2],\, y_1=y_1^{(0)},\, w=u^{(0)}\}$. Since each
pseudopolynomial $P_j$ is a product of irreducible pseudopolynomials of  form
(\ref{pol1}) with conditions (\ref{co}) (see, for example, \cite{H}), we may assume
that all pseudopolynomials $P_j$ in (\ref{rep}) are irreducible. Also, we can rewrite
(\ref{rep}) in the form
 \begin{equation}\label{a}
F(z_1,i'y^{(0)},w)=(z_1-b_1(w))^{t_1}\cdots (z_1-b_k(w))^{t_k} G(z_1,i'y^{(0)},w),
 \end{equation}
where $b_n(w),\,n=1,\dots,k$, are analytic multifunctions in some neighborhood $U$ of
the point $u^{(0)}$.

Since the  functions $\{F(s+iy_1^{(0)},i'y^{(0)},u)\}_{u\in\R^N}$
satisfy the condition of Lemma \ref{2}, we see that the number of
zeros $t_1+t_2+\dots+t_k$ is bounded from above uniformly in
$u\in\R^N\setminus M$. Since an increment of the amplitude of any
linear multiplier along any segment is at most  $\pi$, we see that
the functions $I^+(u)$ and $I^-(u)$ are bounded uniformly in
$u\in\R^N$.

Furthermore, note that  the discriminant $d_P(w)$ of a
pseudopolynomial $P$ of the form (\ref{pol1}) is a holomorphic
function in $U$.  If $P$ is irreducible, then $d_P(w)\not\equiv0$
(see, for example, \cite{H}). Set
 $$
 M_1=\{w\in U:\,F(-1/2,i'y^{(0)},w)
 F(1/2,i'y^{(0)},w)d_{P_1}(w)\cdots d_{P_r}(w)=0\}.
 $$
By the Proposition, $m_N(M_1\cap\R^N)=0$. Take a point $u^{(1)}\in (U\setminus
M_1)\cap\R^N$. Then there is a neighborhood $U_1\subset\C^N$ of $u^{(1)}$ such that
for each point $w\in U_1$ every pseudopolynomial $P_m(z,i'y^{(0)},w)$ has only simple
zeros in $z_1\in\C$. Hence for $w\in U_1$ representation (\ref{a}) holds with  mutual
different functions $b_n(w)$.  We shall prove that each function
$\D_{-1/2<x<1/2}\arg^\pm(x+iy-b_n(u))$ is continuous at points of the set
$U_1\cap\R^N$.

Note that $F(\pm 1/2+iy_1^{(0)},i'y^{(0)},u^{(1)})\neq 0$. Hence
for sufficiently small $U_1$ the functions $b_n(w)$ take $U_1$ to
the set $\{z_1:\,|x_1|<1/2,\,|y_1-y_1^{(0)}|<\delta\}$. Set
 $$
 b'_n(w)= (b_n(w)+\overline{ b_n(\bar w)})/2, \quad
b''_n(w) =(b_n(w)-\overline{ b_n(\bar w)})/2i,
 $$
Note that for all $u\in\R^N$ we have $ b'_n(u) =\Re b_n(u),\,
b''_n(u)=\Im b_n(u)$.

  If $b''_n(w)\not\equiv y_1^{(0)}$ for $w\in U_1$, then, by the
Proposition, $m_N(\{w:\, b''_n(w)=y_1^{(0)}\}\cap\R^N)=0$. Hence
for almost all points $u \in U_1\cap\R^N$ the function
$x_1+iy_1^{(0)}-b_n(u)$ does not vanish for $x_1\in[-1/2,\,1/2]$,
then the functions $\D_{-1/2<x_1<1/2}\arg^+
(x_1+iy_1^{(0)}-b_n(u))$ and
$\D_{-1/2<x_1<1/2}\arg^-(x_1+iy_1^{(0)}-b_n(u))$ are continuous
and coincide almost everywhere on the $U_1\cap\R^N$.

Now consider the case $b''_n(w)\equiv y_1^{(0)}$ for $w\in U_1$.
Since the function $b'_n(u)$ is continuous in $u$, we see that the
function $x_1+iy_1^{(0)}-b_n(u)=x_1-b'_n(u)$ of the variable $x_1$
has exactly one simple zero for all $u$ in a neighborhood $U_2$ of
a point $u^{(2)}\in U_1\cap\R^N$. Hence, $\D_{-1/2<x_1<1/2}\arg^+
(x_1+iy_1-b_n(u))\equiv -\pi$ and
$\D_{-1/2<x_1<1/2}\arg^-(x_1+iy_1-b_n(u))\equiv\pi$ for all $u\in
U_2$.

Since $u^{(0)}$ is an arbitrary point of $[0,\,2\pi]^N\setminus
M$, we see that the functions $I^+(u)$ and $I^-(u)$ are bounded
and continuous almost everywhere in this cube. Therefore, these
functions are integrable in the sense of Riemann over
$[0,\,2\pi]^N$.

Furthermore, by (\ref{shift}), we get
 $$
F(s+iy_1^{(0)},i'y^{(0)},\<\mu^1,x\>,\dots,\<\mu^N,x\>)=P((s,'0)+x+iy^{(0)}).
 $$
 Applying Lemma \ref{1} to the functions $\D_{-1/2<s<1/2}
 \arg^\pm F(s+iy_1^{(0)},i'y^{(0)},u)$,
we obtain that the limits
 \begin{equation}\label{lim1}
\lim_{\min_j(\b_j-\a_j)\to\infty}
{\int_{\Pi^{(p)}(\a,\b)}\D_{-1/2<s<1/2}\arg^+
P((s,'0)+x+iy^{(0)})\,dx\over(\b_1-\a_1)(\b_2-\a_2)\cdots(\b_p-\a_p)}
 \end{equation}
 and
\begin{equation}\label{lim2}
\lim_{\min_j(\b_j-\a_j)\to\infty}
{\int_{\Pi^{(p)}(\a,\b)}\D_{-1/2<s<1/2}\arg^-
P((s,'0)+x+iy^{(0)})\,dx\over(\b_1-\a_1)(\b_2-\a_2)\cdots(\b_p-\a_p)}
 \end{equation}
exist. Since the increments $\D_{-1/2<s<\g}\arg^\pm
P((s,'0)+x+iy^{(0)})=\D_{-1/2<s<\g}\arg^\pm
F((s,'0)+iy^{(0)},\<\mu^1,x\>,\dots,\<\mu^N,x\>)$ are uniformly
bounded  in $s\in\R$, $x\in\R^p$, and $-1/2<\g\le1/2$, we see that
the differences
 $$
\arg^\pm P((s,'0)+x+iy^{(0)})-\int_{ s-1/2}^{s+1/2} \arg^\pm
P(x+iy^{(0)})\,dx_1
 $$
 are uniformly bounded as well. Therefore, up to a uniformly bounded term,
 $$
\arg^\pm P((\b_1,'x)+iy^{(0)})-\arg^\pm
P((\a_1,'x)+iy^{(0)})
 $$
 $$
 =\int^{\b_1+1/2}_{\b_1-1/2}\arg^\pm P((x_1,'x)+iy^{(0)})\,dx_1-
 \int^{\a_1+1/2}_{\a_1-1/2}\arg^\pm
 P((x_1,'x)+iy^{(0)})\,dx_1
 $$
 $$
 =\int^{\b_1}_{\a_1}\D_{-1/2<s<1/2}\arg^\pm P((s,'0)+x+iy^{(0)})\,dx_1.
 $$
 Cosequently,  (\ref{m1}) and
(\ref{m2}) exist and equal (\ref{lim1}) and (\ref{lim2}),
respectively. Theorem is proved.

\bigskip

Mathematical School, Kharkov national university, Swobody sq.4,
Kharkov, 61077 Ukraine

e-mail: Sergey.Ju.Favorov@univer.kharkov.ua

\end{document}